\documentclass[12pt]{amsart}
\usepackage{amstext}
\usepackage{amsmath}
\usepackage{amsfonts, amssymb}
\usepackage{color}
\usepackage{epsfig}
\newtheorem{theorem}{Theorem}
\newtheorem{proposition}[theorem]{Proposition}

\theoremstyle{definition}
\newtheorem{remark}{Remark}

\begin{document}
\title[Stable constant mean curvature hypersurfaces]{Stable constant mean curvature hypersurfaces are area minimizing in small $L^1$ neighborhoods}         
\author[F. Morgan and A. Ros]{Frank Morgan and Antonio Ros}      
\thanks{The first author is partially supported by a National Science Foundation grant and 
 the second author by MEC-FEDER MTM2007-61775 and J. Andaluc\'\i a P06-FQM-01642 grants}         
\maketitle

{\small {\sc Abstract.-}  We prove that a strictly stable constant-mean-curvature hypersurface in a smooth manifold of dimension less than or equal to $7$ is uniquely homologically area minimizing for fixed volume in a small $L^1$ neighborhood.  }

\section{Introduction}

By work of White \cite{W} and Grosse-Brauckman \cite{G}, a strictly stable constant-mean-curvature surface $S_0$ is minimizing in a small neighborhood $U$ of $S_0$ among competitor hypersurfaces $S\subset U$ enclosing the same volume. Assuming $M$ compact, we extend their results to a small $L^1$ neighborhood of $S_0$, i.e., to hypersurfaces $S$ such that $S-S_0$ bounds a region with net volume $0$ and small total volume. 

Stable constant-mean-curvature hypersurfaces in $M$ appear in particular as solutions of the isoperimetric problem; see for instance \cite{Ros}.  
In the case that the ambient space is a flat $3$-torus $T^3$ there is a connection between the isoperimetric problem and the study of mesoscale phase separation phenomena; see Choksi and Sternberg \cite{Cho}. One simple model minimizes the Cahn-Hilliard free energy
\[
\int_{T^3} \left(\frac{\varepsilon^2}{2}|\nabla u|^2 + W(u)\right)dx,
\]
where $W$ is nonnegative with $W(\pm 1)=0$, $u$ represents the concentration of one of the two phases and $\int_{T^3}u\ dx$ is fixed. The local minima of this energy converge as $\varepsilon\rightarrow 0$ to the sharp interface limit given by stable periodic constant-mean-curvature surfaces. This depends on results in $\Gamma$-convergence.

In flat $3$-tori there are some beautiful minimal surfaces, the Schwarz $P$ and $D$ surfaces and the Gyroid $G$ of A. Schoen, which are closely related to complex phases appearing in periodic phase separation. Ross \cite{Ross} has proved that these surfaces are stable for fixed volume and there is a particular interest in providing a mathematical treatment of these complex phases by minimizing locally the Cahn-Hilliard energy or other more sophisticated models. Our $L^1$-minimizing result in this paper gives the necessary tool for such treatment via $\Gamma$-convergence. A different point of view has been considered by Pacard and Ritor\'e \cite{Pacard}.

\vspace{.2cm}

For background in geometric measure theory see Giusti \cite{Gi} and Morgan \cite{M1}.

We are grateful to Rustum Choksi and Peter Sternberg for calling our attention to this problem and to Robert Kohn for helpful conversations.

\section{The proof}

If $S$ and $S'$ are two closed hypersurfaces in $M^{n}$ enclosing regions $\Omega$ and $\Omega'$ respectively, the $L^1$ distance between them is defined as the volume of their symmetric difference: 
\[
||S-S'||_{L^1} =V\left((\Omega-\Omega')\cup(\Omega'-\Omega)\right)
\]   
Note that in order to the define the $L^1$-distance it is enough that $S$ and $S'$ are homologically equivalent, that is, $\partial S=\partial S'$ and $S - S'$ bounds. It is not necessary they are boundaries. 
\begin{figure}[htb]
 \begin{center}
 \includegraphics[width=8.5cm]{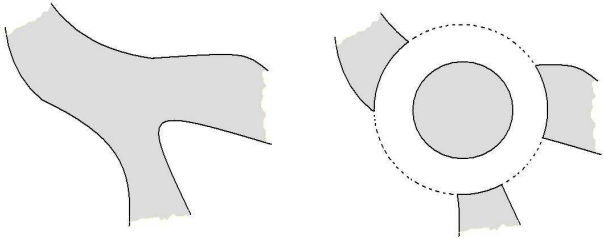}
 \end{center}
 \caption{Proof of the area growth estimate: $S$ at the left. At the right we have a competitor hypersurface which differs from $S$ inside a ball of radius $r$ and encloses the same volume as $S$ in this ball.}
\end{figure}

We will need the following isoperimetric version of the classical result after Fleming \cite[Sect. 5]{F} that for $n\leq 7$, area-minimizing hypersurfaces in the $\mathbb{R}^n$ are hyperplanes.
For $n=3$ da Silveira (\cite{dS}, see also \cite{LR}) proved the result under the weaker hypothesis that $S$ be stable for fixed volume.

\begin{proposition}
\label{cono}
Let $S$ be a hypersurface without boundary in $\mathbb{R}^{n}$, $n\leq 7$, area-minimizing for fixed volume under changes of compact support. Then $S$ is either a 
hypersphere or a hyperplane.
\end{proposition}

\begin{proof} 
If $S$ is compact, $S$ is a hypersphere by the standard isoperimetric inequality. 
Assuming $S$ is not compact, the hypothesis of the proposition implies that $S$ has constant mean curvature and is stable for fixed volume. Given $r>0$, inside the ball $B(r)$ about a fixed point of $S$, replace the region bounded by $S$ by a ball $B(\rho)$, $0\leq\rho\leq r$, of the same volume as in 
Figure 1. The resulting area inside $B(r)$ is at 
most twice the area of the hypersphere $\partial B(r)$. Hence the original area inside the ball 
of the minimizer $S$ is at most $C r^{n-1}$, for some $C$. 
By Cheung \cite{C} $S$ has mean curvature $0$. By monotonicity of the 
mass ratio \cite[Cor. 5.1(3) p. 446]{A}, the area
divided by $\alpha_n r^{n-1}$, where $\alpha_n$ is the 
volume of the unit ball in $\mathbb{R}^{n-1}$, is nondecreasing 
in $r$, varying from $1$ as $r$ approaches $0$ to a limit 
$C_0$, as $r$ approaches infinity. Therefore 
homothetic contractions, restricted to balls 
about the origin, have area bounded below and 
above, so that by compactness \cite[5.5 and remark 
page 88]{M1}, a subsequence converges to a nonzero 
limit, which has constant area ratio $C_0$ and is 
therefore a cone \cite[Cor. 5.1(2)]{A}. Since the cone 
minimizes area for given volume and $n\leq 7$, by 
regularity \cite{M2} the cone must be a hyperplane 
(with multiplicity $1$ because it is the boundary 
of a region) and $C_0 = 1$. Hence likewise $S$ has 
constant mass ratio $1$ and must be a hyperplane.
\end{proof}

Now we prove our main result.

\begin{theorem} In a smooth closed Riemannian 
manifold of dimension $n\leq 7$, let $S_0$ be a smooth 
constant-mean-curvature hypersurface, possibly with boundary, with 
positive second variation for fixed volume and boundary. Then 
$S_0$ is uniquely homologically area minimizing for fixed volume among hypersurfaces in a small $L^1$ neighborhood. 
\end{theorem}

In particular, if $S_0$ bounds a region, then it minimizes area among hypersurfaces $S$ enclosing the same volume with $||S-S_0||_{L^1}$ small.
It is not necessary to assume that $S_0$ is a boundary. Our proof gives that $S_0$ minimizes among competitors $S$ such that $\partial S=\partial S_0$ and 
$S - S_0$ bounds net oriented volume $0$.

\begin{proof}
Denote area, volume, and mean curvature by 
$A$, $V$, and $H$. The subscript $0$ refers to $S_0$. Our hypersurface $S_0$ has positive second variation 
under smooth variations which fix volume (or 
equivalently under smooth variations which fix 
volume to first order). By Grosse-Brauckmann \cite[Lemma 5]{G}, for some $C > 0$, $S_0$ has positive smooth 
second variation for the energy
\[
        F = A + H_0 V + (C/2)(V - V_0)^2
\]
under general smooth variations. As 
Grosse-Brauckman \cite[last paragraph]{G}
%[last paragraph]{G} 
points out, 
White \cite[Thm. 3]{W}  applies to show that $S_0$ uniquely minimizes $F$ in a neighborhood. To see this, let $\omega$ 
be a smooth differential form which over 
homologous surfaces gives the volume enclosed 
with $S_0$, such that $C\omega$ is small in a 
neighborhood of $S_0$ \cite[end of Intro.]{W}. To apply 
\cite[Thm. 3]{W}, take $F$ to be the area integrand, $F_1 
= F + C \omega$, $F_2 = F$, and $\phi(x,y) = 
(x-y)^2/2C$. By \cite[Thm. 3]{W}, $S_0$ uniquely 
minimizes $F$ in a small neighborhood $U$ of its 
support. In particular, for fixed volume, $S_0$ 
uniquely minimizes $A$ in $U$.

  To obtain a contradiction, suppose that 
there is a sequence of surfaces $S_i$ of less area 
than $S_0$ converging in $L^1$ to 
$S_0$ and enclosing net signed volume $0$ with $S_0$. We may assume that 
$S_i$ minimizes area for fixed  $||S_i - S_0||_{L^1} = 
\varepsilon_i \rightarrow 0$. On the exterior and on the 
interior of $S_0$, $S_i$ minimizes area for fixed 
volume; therefore $S_i$ is a smooth constant mean curvature
surface \cite[Cor. 3.7]{M2}  (although the exterior 
constant need not equal the interior constant; we 
assert no regularity at points of $S_0$). By the 
first paragraph of this proof, each $S_i$ strays 
outside $U$. By replacing $S_i$ by a subsequence, we 
may assume that each $S_i$ strays outside of $U$ 
always on the same side of $S_0$ or on both sides 
of $S_0$.

        Hence by monotonicity, on a relevant side 
of $S_0$, the curvature of the sequence $S_i$ is not bounded in $M-U$. Indeed, if 
the mean curvature were bounded, then by 
monotonicity of the mass ratio \cite[Cor. 5.1(3) 
p. 446 and Rmk. 4.4]{A} , the area of $S_i$ outside a 
smaller neighborhood $U'$ would be bounded below by 
some positive constant $\delta$, and then
\[
        A(S_0) \leq \liminf A(S_i) - \delta \leq A(S_0) - \delta,
\]
the desired contradiction.

        Choose a point outside of $U$ on a relevant 
side of $S_0$ on each $S_i$ of maximum $|II|^2$  (the 
sum of the squares of the principal curvatures) 
and blow up the picture to make $|II|^2 = 1$. A 
limit is minimizing for fixed volume in $\mathbf{R}^{n}$ and 
hence must be a round sphere by Proposition \ref{cono}. 
Hence for some large $i$, $S_i$ includes a small, 
nearly round sphere partly outside $U$. We may 
assume that there are no other points of $S_i$ 
outside $U$ on that side of $S_0$, since otherwise we 
could repeat the argument on $S_i$ minus the first 
sphere and obtain a second such sphere, while 
combining them as one sphere would do better. 
Hence on each side of $S_0$, there is at most one 
such sphere partly outside $U$. For a constant $c_n$ 
depending only on the dimension $n$, the total area 
and volume of such spheres satisfy $a > c_n v^{(n-1)/{n}}$.

        Let $T_i$ be $S_i$ minus such spheres, so 
that $T_i$ lies in 
%an $\eta$ 
the neighborhood $U$ of $S_0$. Now
\[
                F(T_i) < A(S_i) - c_n v^{(n-1)/{n}} + |H_0| v + (C/2)v^2 < A(S_i)
\]
for small $v$ and hence for large $i$. Then
\[
        F(T_i) < A(S_i) < A(S_0) = F(S_0),
\]
a contradiction of the fact that $S_0$ minimizes $F$ in $U$.
\end{proof}

\begin{remark} 
For minimal surfaces, the result also holds without volume constraints. The same proof holds, with simplifications.
\end{remark}

\begin{remark}
In applications, as the Cahn-Hilliard problem in flat $3$-tori, it is important to consider the case where the ambient space $M$ has non-trivial isometry group and $S_0$ is a closed constant-mean-curvature hypersurface with positive second variation orthogonal to the isometries, for fixed volume. Our proof applies in this case without changes as White \cite{W} also holds. 
White's proof observes that a sequence of other minimizers in shrinking physical neighborhoods of $S$ are {\it almost minimizing} and hence H\"older differentiable manifolds that converge H\"older differentiably to $S$, contradicting the positive second variation of $S$. 
In the presence of isometries, one may translate the nearby minimizers to be graphs of functions orthogonal to the isometries to obtain the same contradiction. 

In particular, as Ross \cite{Ross} proved that the $P$, $D$ and $G$ minimal surfaces have positive second variation for any direction orthogonal to the ones induced by the isometries of $T^3$, if follows that they uniquely minimize area, up to isometries, for fixed volume in a small $L^1$ neighborhood in the ambient $3$-torus. For a description of stable constant-mean-curvature surfaces with fixed volume in flat $3$-tori see Ros \cite{Ros1}.
\end{remark}

{\footnotesize
\noindent
{\sc Frank Morgan}\\
{\tt Frank.Morgan@williams.edu}\\ 
{ Department of Mathematics  and Statistics}\\
{ Williams College}\\
{ Williamstown, Massachusetts 01267 USA\\
}}

{\footnotesize
\noindent
{\sc Antonio Ros}\\
{\tt aros@ugr.es}\\ 
{ Departamento de Geometr\'\i a y Topolog\'\i a}\\
{ Facultad de Ciencias - Universidad de Granada}\\
{ 18071  Granada  Spain}}

\end{document}